\documentclass{article}
\usepackage{amssymb}
\usepackage{amscd}
\usepackage{amsfonts}
\usepackage{amsmath}
\usepackage{graphicx}
\usepackage{cite}
\usepackage{enumerate}
\usepackage{subfigure}
\usepackage[latin5]{inputenc}
\usepackage{authblk}
\usepackage{float}
\usepackage{graphicx}
\usepackage{epstopdf}
\usepackage{amscd,amsmath,amssymb,amsthm,amsfonts,epsfig,graphics}
\DeclareMathOperator{\sech}{sech}
\usepackage[english]{babel}
\usepackage{amsthm}
\newtheorem*{theorem}{Theorem}%[section]
\begin{document}
\author[1]{Ozlem Ersoy\thanks{Ozlem Ersoy, Telephone: +90 222 2393750, Fax: +90 222 2393578, E-mail:ozersoy@ogu.edu.tr}$^,$}
\author[1]{Idiris Dağ}
\author[2]{Alper Korkmaz}
\affil[1]{Department of Mathematics-Computer, Eskişehir Osmangazi University, 26480, Eskişehir, Turkey.}
\affil[2]{Department of Mathematics, Çankırı Karatekin University, 18200, Çankırı, Turkey. }
\title{Solitary wave simulations of the Boussinesq Systems}
\maketitle
\begin{abstract}
In the study, the collocation method based on exponential cubic B-spline functions is proposed to solve one dimensional Boussinesq systems numerically. Two initial boundary value problems for Regularized and Classical Boussinesq systems modeling motion of traveling waves are considered. The accuracy of the method is validated by measuring the error between the numerical and analytical solutions. The numerical solutions obtained by various values of free parameter $p$ are compared with some solutions in literature. 
\end{abstract}
\textit{Keywords:} Boussinesq sytems; solitary waves; exponential cubic B-spline; collocation.

\noindent
\textbf{MSC2010}:\, \textit{35Q99;35C07;76B25;65D07;65M70.}

\section{Boussinesq Systems}

In the light of the d' Alembert solution, describing two distinct waves moving in the opposite directions for the Cauchy problem for one dimensional wave equation, many physical problems modeled by linear and nonlinear partial differential equations have been solved in wave forms covering traveling waves, solitary waves, harmonic waves, etc. Waves are an important solution class for the model problems in physics, chemistry, and many fields of engineering. This solution type is widely constructed for different models in various media covering solids and fluids. Class of surface bell-shaped solitary waves, having long amplitude when compared with its width, is one of the most prominent classes of waves. In addition to linear equations, solitary wave solutions have also been studied for well known nonlinear equations and systems such as Korteweg-deVries(KdV), Schrödinger, Boussinesq, etc. Having analytical solutions only for some particular cases including additional restrictions and conditions, many numerical methods have been developed to obtain solutions for nonlinear problems. 

\noindent
Dougalis et al.\cite{dougalis} studied the numerical behavior of solitary waves of a member of Boussinesq systems family (Bona-Smith system). In that study, long time stability and possible blow-up solutions were examined under small and large perturbations covering perturbation of amplitudes, system coefficients, etc. Numerical models were constructed with the Galerkin-finite element method on the space \textit{$S_h$} of smooth, periodic, cubic splines.

\noindent
Quintic B-spline collocation tecnique based on Crank-Nicolson formulation was set up for numerical solutions of Boussines type coupled-BBM system\cite{behzadi}. Two initial boundary value problems describing motion of single solitary wave and interaction of solitary waves were simulated by the proposed method. The numerical results were compared with the analytical ones.

\noindent 
The Galerkin-finite element periodic B-spline method was used to approximate the cnoidal and solitary wave solutions of periodic initial value problem for four variable Boussinesq system\cite{anton1}. 

\noindent
In the study\cite{karasozen}, variations of Boussinesq system, namely coupled Korteweg-de Vries (KdV) systems in Hamiltonian form, were integrated using the energy preserving average vector field scheme. The travelling waves for the KdV-KdV systems were also studied numerically by Bona et al. using unconditionally stable periodic splines Galerkin method with Gauss-Legendre implicit Runge-Kutta method of stage two\cite{bona2010}.  

\noindent
Antonopoulos et al. simulated the numerical solutions of three initial boundary value problems one parameter Bona-Smith systems\cite{anton2}. The problems with homogenous Dirichlet, reflection, and periodic boundary conditions were integrated by Galerkin-finite element combined with fourth order explicit Runge-Kutta method. 

\noindent
A meshless radial basis collocation algorithm was set up to simulate nonlinear dispersive waves propagating in two ways of Boussinesq system by Su\'{a}rez and Morales\cite{suarez}. The solutions for three initial boundary value problems modeling motion of single solitary and interaction of solitary waves were demonstrated with the proposed method.

\noindent
Aksoy presented Taylor-Collocation extended cubic B-spline methods for the numerical solutions of the Boussinesq systems\cite{aksoy}. The initial boundary value problems modeling motion of single solitary waves were investigated by using various free parameters in extended cubic B-spline functions.
   
\noindent
In this study, the authors investigate the numerical solutions for two particular initial boundary value problems for the Classical Boussinesq System (CBS) and the Regularized Boussinesq System (RBS). A collocation method based on exponential cubic B-spline functions is constructed to obtain the numerical solutions for both systems. The details of both problems and numerical methodology will be explained in the following sections.

\noindent
Consider the Boussinesq system modeling bi-directional propagation of nonlinear dispersive long surface waves with small amplitudes in a channel of the form
\begin{equation} \label{eq:bs}
\begin{aligned}
V_{t}+U_{x}+(VU)_{x}+s_{1}U_{xxx}-s_{2}V_{xxt}&=0 \\ 
U_{t}+V_{x}+UU_{x}+s_{3}V_{xxx}-s_{4}U_{xxt}&=0%
\end{aligned}%
\end{equation}%
where $x$ denotes the distance along the channel scaled by still water of depth $h_0$ and $t$ is the time scaled by $\sqrt{\frac{h_0}{g}}$, with the gravitational acceleration $g$\cite{bona97,Chena}. In order to equalize the order of the effects of both nonlinearity and dispersion, the Stokes number $S=s_1\lambda^2 / h^3_0$ is assumed of order one with the small maximum deviation of surface $s_1$ relative to $h_0$ and $\lambda>h_0$. In the model (\ref{eq:bs}), $V(x,t)$ is the scaled, dimensionless deviation of the water surface, $U(x,t)$ is the scaled, dimensionless horizontal velocity of height $\theta h_0$, $\theta \in [0,1]$ above the bottom of the channel\cite{Chena}.
In the system (\ref{eq:bs}), the real constants $s_1,s_2,s_3$ and $s_4$ are defined by;
\begin{equation}
\begin{aligned}
s_{1}&=\frac{1}{2}(\theta ^{2}-\frac{1}{3})\lambda \\ 
s_{2}&=\frac{1}{2}(\theta ^{2}-\frac{1}{3})(1-\lambda ) \\ 
s_{3}&=\frac{1}{2}(1-\theta ^{2})\mu \\ 
s_{4}&=\frac{1}{2}(1-\theta ^{2})(1-\mu )%
\end{aligned}
\label{BST3}
\end{equation}
where $\lambda$ and $\mu$ are real\cite{bona97,BonaChen,Bona2002}.
\\
Particular choice of the parameters $\theta^2=\frac{1}{3}, \mu=0$ and arbitrary $\lambda$ generates the CBS from the system (\ref{eq:bs})\cite{BS,BonaChen,Chena} as;
\begin{equation}
\begin{aligned}
V_{t}+U_{x}+(VU)_{x}=0 \\ 
U_{t}+V_{x}+UU_{x}-\frac{1}{3}U_{xxt}=0%
\end{aligned}
\label{eq:CBS}
\end{equation}
When the parameters are chosen as $\theta^2=\frac{2}{3}, \lambda=0, \mu=0$, the Boussinesq system (\ref{eq:bs}) reduces to the RBS of the form \cite{Bona2002};
\begin{equation}
\begin{aligned}
V_{t}+U_{x}+(VU)_{x}-\frac{1}{6}V_{xxt}&=0 \\ 
U_{t}+V_{x}+UU_{x}-\frac{1}{6}U_{xxt}&=0%
\end{aligned}
\label{eq:RBS}
\end{equation}
\begin{theorem}
Solitary wave solutions for the Boussinesq system (\ref{eq:bs}) is given as;
\begin{equation}
\begin{aligned}
V(x,t)&=V_{0}\sech ^{2}(\lambda (x+x_{0}-C_{s}t)) \\ 
U(x,t)&=\pm \sqrt{\dfrac{3}{V_{0}+3}}V_{0}\sech ^2{(\lambda (x+x_{0}-C_{s}t))}%
\end{aligned}
\label{eq:TWS}
\end{equation}
where
\begin{equation*}
C_{s}=\frac{3+2V_{0}}{\pm \sqrt{3(3+V_{0})}},\text{ }\lambda =\frac{1}{2}%
\sqrt{\frac{2V_{0}}{3(s_{1}-s_{2})+2s_{2}(V_{0}+3)}}
\end{equation*}%
and $V_{0}$ can be any constant satisfying,
\begin{enumerate}[i.]
\item $V_{0}=\dfrac{3(1-2\kappa)}{2\kappa}$ when $s_{1}-s_{2}+2s_{4}\neq 0,\text{ }\kappa=\dfrac{-s_{2}+s_{3}+2s_{4}}{s_{1}-s_{2}+2s_{4}}>0\text{ and }(\kappa-1/2)\left[ (s_{2}-s_{1})\kappa-s_{2}\right] >0$
\item $0<V_{0}<\infty$ when $s_{1}=s_{2}=s_{3}>0,\text{ }s_{4}=0$
\item $-3\leq V_{0}<0$ when $s_{1}=s_{2}=s_{3}<0,\text{ }s_{4}=0$
\item $V_{0}>-3\text{ and }\dfrac{3}{V_{0}+3}\notin \left[1,\dfrac{s_{2}}{s_{4}}\right]$ when $s_{1}-s_{2}+s_{4}=0,\text{ }s_{1}=s_{3},\text{ }s_{4}>0$
\item $V_{0}>-3\text{ and }\dfrac{3}{V_{0}+3}\in \left[ 1,\dfrac{s_{2}}{s_{4}}\right]$ when $s_{1}-s_{2}+s_{4}=0,\text{ }s_{1}=s_{3},\text{ }s_{4}<0$
\end{enumerate}
It should be noted that the constants $V_{0}$ and $\pm \sqrt{\frac{3}{V_{0}+3}}V_{0}$ denote the amplitudes of solitary waves and the constant $C_{s}$ denotes
the velocity of solitary waves in the analytical solution given above\cite{Chenb}.
\end{theorem}
\noindent
Numerical solutions will be obtained for both initial boundary value problems with the homogenous Neumann boundary conditions at both ends
\begin{equation}
\begin{aligned}
\left . \dfrac{\partial V(x,t)}{\partial x}\right |_{x=a}&=0, \, &\left . \dfrac{\partial V(x,t)}{\partial x}\right |_{x=b}=0, \, t\geq 0 \\ 
\left . \dfrac{\partial U(x,t)}{\partial x}\right |_{x=a}&=0, \, &\left . \dfrac{\partial U(x,t)}{\partial x}\right |_{x=b}=0, \, t\geq 0%
\end{aligned}
\label{BSTsin}
\end{equation}
over the finite problem interval $[a,b]$ are combined with both the CBS and the RBS. The initial conditions for both problems will be chosen in the following sections.

\section{Exponential Cubic B-spline Collocation Method(ECC)}
Let $\pi$ be a uniformly distributed grids of the finite interval $[a,b]$, such as,
$$\pi : x_m=a+mh, m=0,1, \ldots N$$
where $h=\frac{b-a}{N}$. Then, the exponential cubic B-splines are defined as;
\begin{equation}
B_{m}(x)=\left\{ 
\begin{array}{lcc}
b_{2}\left ( (x_{m-2}-x)-\frac{1}{\zeta} \sinh(\zeta(x_{m-2}-x)) \right ) & , & [x_{m-2},x_{m-1}] \\ 
a_{1}+b_{1}(x_{m}-x)+c_{1}\exp(\zeta(x_{m}-x))+d_{1}\exp(-\zeta(x_{m}-x)) & , & [x_{m-1},x_{m}] \\ 
a_{1}+b_{1}(x-x_{m})+c_{1}\exp(\zeta(x-x_{m}))+d_{1}\exp(-\zeta(x-x_{m})) & , & [x_{m},x_{m+1}] \\ 
b_{2}\left ( (x-x_{m+2})-\frac{1}{\zeta} \sinh(\zeta(x-x_{m+2})) \right ) & , & [x_{m+1},x_{m+2}] \\ 
0 & , & otherwise%
\end{array}%
\right.  \label{ecbs}
\end{equation}%
where 
$$
\begin{array}{l}
a_{1}=\dfrac{\zeta h\cosh(\zeta h)}{\zeta h\cosh(\zeta h)-\sinh(\zeta h)}, \\
b_{1}=\dfrac{\zeta }{2} \dfrac{\cosh(\zeta h)(\cosh(\zeta h)-1)+\sinh^2(\zeta h)}{(\zeta h\cosh(\zeta h)-\sinh(\zeta h))(1-\cosh(\zeta h))}, \\
b_{2}=\dfrac{\zeta }{2(\zeta h\cosh(\zeta h)-\sinh(\zeta h))}, \\
c_{1}=\dfrac{1}{4} \dfrac{\exp(-\zeta h)(1-\cosh(\zeta h))+\sinh(\zeta h)(\exp(-\zeta h)-1))}{(\zeta h\cosh(\zeta h)-\sinh(\zeta h))(1-\cosh(\zeta h))}, \\
d_{1}=\dfrac{1}{4} \dfrac{\exp(\zeta h)(\cosh(\zeta h)-1)+\sinh(\zeta h)(\exp(\zeta h)-1))}{(\zeta h\cosh(\zeta h)-\sinh(\zeta h))(1-\cosh(\zeta h))},
\end{array}
$$ 
where $\zeta $ is a real parameter\cite{mccartin}. Each exponential cubic B-spline $B_m(x)$ has two continuous principle derivatives defined in the interval $[x_{m-1},x_{m+2}]$ and $B_m(x)$ itself and its two principle derivatives vanish out of this interval. The set $\{B_{-1}(x), B_{0}(x), ..., B_{N+1}(x)\}$ constitutes a basis for the functions defined over the interval $[a,b]$. Since introduced by McCartin, exponential cubic B-spline functions have been used to solve some engineering and physics problems numerically \cite{moh1,moh2,ozkdv}. 
\section{Numerical Approximation}
Let $U(x,t)$ and $V(x,t)$ be two approximate solutions defined by;
\begin{subequations}
\begin{equation} \label{app1}
U(x,t) \cong \sum_{m=-1}^{N+1}\delta _{m}B_{m}(x)
\end{equation}%
\begin{equation} \label{app2}
V(x,t) \cong \sum_{m=-1}^{N+1}\phi _{m}B_{m}(x) 
\end{equation}%
\end{subequations}
where $\delta _{m}$ and $\phi _{m}$ are parameters dependent on time variable. Those parameters will be determined by the collocation method. 
Then, the two principle derivatives of  ${U}(x,t)$ and ${V}(x,t)$ are;
\begin{subequations}
\begin{equation} \label{u1}
{U}^{'}(x,t)\cong\sum_{m=-1}^{N+1}\delta _{m}B^{'}_{m}(x) 
\end{equation}%
\begin{equation} \label{u2}
{U}^{''}(x,t)\cong\sum_{m=-1}^{N+1}\delta _{m}B^{''}_{m}(x)
\end{equation}%

\begin{equation} \label{v1}
{V}^{'}(x,t)\cong\sum_{m=-1}^{N+1}\phi _{m}B^{'}_{m}(x) 
\end{equation}%
\begin{equation} \label{v2}
{V}^{''}(x,t)\cong\sum_{m=-1}^{N+1}\phi _{m}B^{''}_{m}(x)
\end{equation}
\end{subequations}
Using Eq. (\ref{app1})-(\ref{v2}), the approximate nodal values of the functions ${U}(x,t)$\ and ${V}(x,t)$ and their two principle derivatives are calculated as;
\scriptsize
\begin{equation} \label{app3}
\begin{aligned}
U_{m}=U(x_{m},t)& \cong\dfrac{\sinh{(\zeta h)}-\zeta h}{2(\zeta h\cosh{(\zeta h)}-\sinh{(\zeta h)})}\delta _{m-1}+\delta _{m}+\dfrac{\sinh{(\zeta h)}-\zeta h}{2(\zeta h\cosh{(\zeta h)}-\sinh{(\zeta h)})}\delta _{m+1} \\
U_{m}^{'}=U^{' }(x_{m},t)& \cong\dfrac{\zeta (1-\cosh{(\zeta h)})}{2(\zeta h\cosh{(\zeta h)}-\sinh{(\zeta h)})}\delta _{m-1}+\dfrac{\zeta (\cosh{(\zeta h)}-1)}{2(\zeta h\cosh{(\zeta h)}-\sinh{(\zeta h)})}\delta _{m+1} \\
U_{m}^{''}=U^{'' }(x_{m},t)& \cong\dfrac{\zeta ^{2}\sinh{(\zeta h)}}{2(\zeta h\cosh{(\zeta h)}-\sinh{(\zeta h)})}\delta _{m-1}-\dfrac{\zeta ^{2}\sinh{(\zeta h)}}{\zeta h\cosh{(\zeta h)}-\sinh{(\zeta h)}}\delta _{m}+\dfrac{\zeta ^{2}\sinh{(\zeta h)}}{2(\zeta h\cosh{(\zeta h)}-\sinh{(\zeta h)})}\delta_{m+1}\\
V_{m}=V(x_{m},t)&\cong\dfrac{\sinh{(\zeta h)}-\zeta h}{2(ph\cosh{(\zeta h)}-\sinh{(\zeta h)})}\phi _{m-1}+\phi _{m}+\dfrac{\sinh{(\zeta h)}-\zeta h}{2(\zeta h\cosh{(\zeta h)}-\sinh{(\zeta h)})}\phi _{m+1}\\ 
V_{m}^{'}= V^{\prime }(x_{m},t)&\cong\dfrac{\zeta (1-\cosh{(\zeta h)})}{2(\zeta h\cosh{(\zeta h)}-\sinh{(\zeta h)})}\phi _{m-1}+\dfrac{\zeta (\cosh{(\zeta h)}-1)}{2(\zeta h\cosh{(\zeta h)}-\sinh{(\zeta h)})}\phi _{m+1} \\
V_{m}^{''}=V^{\prime \prime }(x_{m},t)&\cong\dfrac{\zeta ^{2}\sinh{(\zeta h)}}{2(\zeta h\cosh{(\zeta h)}-\sinh{(\zeta h)})}\phi _{m-1}-\dfrac{\zeta ^{2}\sinh{(\zeta h)}}{\zeta h\cosh{(\zeta h)}-\sinh{(\zeta h)}}\phi _{m}+\dfrac{\zeta ^{2}\sinh{(\zeta h)}}{2(\zeta h\cosh{(\zeta h)}-\sinh{(\zeta h)})}\phi_{m+1}
\end{aligned}
\end{equation}
\normalsize
\section{Discretization of the Boussinesq System}
Under the assumption $s_{1}=s_{3}=0$, substituting finite difference approximations into the terms containing the derivatives with respect to the time variable $t$ and Crank-Nicolson approximations instead of the remaining terms reduces the Boussinessq system to;
\begin{equation} \label{nl}
\begin{aligned}
\dfrac{V^{n+1}-V^{n}}{\Delta t}+\dfrac{U_{x}^{n+1}+U_{x}^{n}}{2}+\dfrac{(VU)_{x}^{n+1}+(VU)_{x}^{n}}{2}-s_{2}\dfrac{V_{xx}^{n+1}-V_{xx}^{n}}{\Delta t}&=0 \\ 
\dfrac{U^{n+1}-U^{n}}{\Delta t}+\dfrac{V_{x}^{n+1}+V_{x}^{n}}{2}+\dfrac{(UU_{x})^{n+1}+(UU_{x})^{n}}{2}-s_{4}\dfrac{U_{xx}^{n+1}-U_{xx}^{n}}{\Delta t}&=0
\end{aligned}
\end{equation}
where $U^{n+1}=U(x,t_{n}+\Delta t)$ and $V^{n+1}=V(x,t_{n}+\Delta t)$.
The nonlinear terms $(VU)_{x}^{n+1}$ and $(UU_{x})^{n+1}$ are linearized by using Rubin{\&}Gravis' technique\cite{rubin} as;
\begin{equation}
\begin{aligned}
(VU)_{x}^{n+1} & \cong (V_{x}U)^{n+1}+(VU_{x})^{n+1} \cong V_{x}^{n+1}U^{n}+V_{x}^{n}U^{n+1}-V_{x}^{n}U^{n}+V^{n+1}U_{x}^{n}+V^{n}U_{x}^{n+1}-V^{n}U_{x}^{n} \\ 
(UU_{x})^{n+1} & \cong U^{n+1}U_{x}^{n}+U^{n}U_{x}^{n+1}-U^{n}U_{x}^{n}
\end{aligned}
\label{linearization}
\end{equation}
After linearization of system (\ref{nl}), the approximates (\ref{app1}), (\ref{app2}) and (\ref{app3}) to the functions and derivatives are substituted into the linearized form of the system. The resulting equation system are obtained at each grid point as;
\begin{equation} \label{disc1}
\begin{aligned}
&\nu _{m1}\delta _{m-1}^{n+1}+\nu _{m2}\phi _{m-1}^{n+1}+\nu _{m3}\delta_{m}^{n+1}+\nu _{m4}\phi _{m}^{n+1}+\nu _{m5}\delta _{m+1}^{n+1}+\nu
_{m6}\phi _{m+1}^{n+1} \\
=&\nu _{m7}\delta _{m-1}^{n}+\nu _{m8}\phi _{m-1}^{n}+\nu _{m9}\delta_{m}^{n}+\nu _{m10}\phi _{m}^{n}+\nu _{m7}\delta _{m+1}^{n}+\nu _{m11}\phi
_{m+1}^{n}
\end{aligned}
\end{equation}
and
\begin{equation} \label{disc2}
\begin{aligned}
&\nu _{m12}\delta _{m-1}^{n+1}+\nu _{m13}\phi _{m-1}^{n+1}+\nu _{m14}\delta_{m}^{n+1}+\nu _{m15}\phi _{m}^{n+1}+\nu _{m16}\delta _{m+1}^{n+1}+\nu
_{m17}\phi _{m+1}^{n+1} \\
=&\nu _{m18}\delta _{m-1}^{n}+\nu _{m19}\phi _{m-1}^{n}+\nu _{m20}\delta_{m}^{n}+\nu _{m21}\phi _{m}^{n}+\nu _{m22}\delta _{m+1}^{n}+\nu _{m19}\phi_{m+1}^{n}
\end{aligned}
\end{equation}
where
\begin{equation*}
\begin{array}{lll}
\nu _{m1}=\left( \dfrac{2}{\Delta t}+K_{2}\right) \alpha _{1}+K_{1}\beta
_{1}-\dfrac{2s_{4}}{\Delta t}\gamma _{1}, &  & \nu _{m7}=\dfrac{2}{\Delta t}
\alpha _{1}-\dfrac{2s_{4}}{\Delta t}\gamma _{1}, \\ 
\nu _{m2}=\beta _{1}, &  & \nu _{m8}=-\beta _{1}, \\ 
\nu _{m3}=\left( \dfrac{2}{\Delta t}+K_{2}\right) -\dfrac{2s_{4}}{\Delta t}
\gamma _{2}, &  & \nu _{m9}=\dfrac{2}{\Delta t}-\dfrac{2s_{4}}{\Delta t}
\gamma _{2}, \\ 
\nu _{m4}=0, &  & \nu _{m10}=0, \\ 
\nu _{m5}=\left( \dfrac{2}{\Delta t}+K_{2}\right) \alpha _{1}-K_{1}\beta
_{1}-\dfrac{2s_{4}}{\Delta t}\gamma _{1}, &  & \nu _{m11}=\beta _{1}, \\ 
\nu _{m6}=-\beta _{1}, &  &  \\ 
&  &  \\ 
\nu _{m12}=L_{2}\alpha _{1}+(1+L_{1})\beta _{1}, &  & \nu _{m18}=-\beta _{1},
\\ 
\nu _{m13}=\left( \dfrac{2}{\Delta t}+K_{2}\right) \alpha _{1}+K_{1}\beta
_{1}-\dfrac{2s_{2}}{\Delta t}\gamma _{1}, &  & \nu _{m19}=\dfrac{2}{\Delta t}
\alpha _{1}-\dfrac{2s_{2}}{\Delta t}\gamma _{1}, \\ 
\nu _{m14}=L_{2}, &  & \nu _{m20}=0, \\ 
\nu _{15}=\left( \dfrac{2}{\Delta t}+K_{2}\right) -\dfrac{2s_{2}}{\Delta t}
\gamma _{2}, &  & \nu _{m21}=\dfrac{2}{\Delta t}-\dfrac{2s_{2}}{\Delta t}\gamma _{2}, \\ 
\nu _{m16}=L_{2}\alpha _{1}-(1+L_{1})\beta _{1}, &  & \nu _{m22}=\beta _{1},
\\ 
\nu _{m17}=\left( \dfrac{2}{\Delta t}+K_{2}\right) \alpha _{1}-K_{1}\beta
_{1}-\dfrac{2s_{2}}{\Delta t}\gamma _{1}, &  & 
\end{array}%
\end{equation*}%
with
\begin{equation*}
\begin{aligned}
K_{1}&=\alpha _{1}\delta _{m-1}^{n}+\delta _{m}^{n}+\alpha _{1}\delta_{m+1}^{n} \\
L_{1}&=\alpha _{1}\phi _{m-1}^{n}+\phi _{m}^{n}+\alpha _{1}\phi_{m+1}^{n} \\
K_{2}&=\beta _{1}\delta _{m-1}^{n}-\beta _{1}\delta _{m+1}^{n} \\
L_{2}&=\beta_{1}\phi _{m-1}^{n}-\beta _{1}\phi _{m+1}^{n}\\
\alpha _{1}&=\dfrac{\sinh{(\zeta h)}-\zeta h}{2(\zeta h\cosh{(\zeta h)}-\sinh{(\zeta h)})} \\
\beta _{1}&=\dfrac{\zeta (1-\cosh{(\zeta h)})}{2(\zeta h\cosh{(\zeta h)}-\sinh{(\zeta h)})} \\
\gamma _{1}&=\dfrac{\zeta ^{2}\sinh{(\zeta h)}}{2(\zeta h\cosh{(\zeta h)}-\sinh{(\zeta h)})}\\
\gamma _{2}&=-\dfrac{\zeta ^{2}\sinh{(\zeta h)}}{\zeta h\cosh{(\zeta h)}-\sinh{(\zeta h)}}
\end{aligned}
\end{equation*}
The discretized system (\ref{disc1})-(\ref{disc2}) is written in matrix notation as;

\begin{equation} \label{mat}
	\mathbf{Ax^{n+1}=Bx^n}
\end{equation}
where
\begin{equation*}
\mathbf{A=}%
\begin{bmatrix}
\nu _{m1} & \nu _{m2} & \nu _{m3} & \nu _{m4} & \nu _{m5} & \nu _{m6} &  & 
&  &  \\ 
\nu _{m12} &\nu _{m13} & \nu _{m14} & \nu _{m15} & \nu _{m16} & \nu _{m17} 
&  &  &  &  \\ 
&  & \nu _{m1} & \nu _{m2} & \nu _{m3} & \nu _{m4} & \nu _{m5} & \nu _{m6} & 
&  \\ 
&  & \nu _{m12} & \nu _{m13} & \nu _{m14} & \nu _{m15} & \nu _{m16} & \nu _{m17} &  &  \\ 
&  &  & \ddots  & \ddots  & \ddots  & \ddots  & \ddots  & \ddots  &  \\ 
&  &  &  & \nu _{m1} & \nu _{m2} & \nu _{m3} & \nu _{m4} & \nu _{m5} & \nu
_{m6} \\ 
&  &  &  & \nu _{m12} & \nu _{m13} & \nu _{m14} & \nu _{m15} & \nu _{m16} &\nu _{m17}
\end{bmatrix}%
\end{equation*}%
,
\begin{equation*}
\mathbf{B=}%
\begin{bmatrix}
\nu _{m7} & \nu _{m8} & \nu _{m9} & \nu _{m10} & \nu _{m7} & \nu _{m11} &  & 
&  &  \\ 
\nu _{m18} & \nu _{m19} & \nu _{m20} & \nu _{m21} & \nu _{m22} & \nu _{m19}
&  &  &  &  \\ 
&  & \nu _{m7} & \nu _{m8} & \nu _{m9} & \nu _{m10} & \nu _{m7} & \nu _{m11}
&  &  \\ 
&  & \nu _{m18} & \nu _{m19} & \nu _{m20} & \nu _{m21} & \nu _{m22} & \nu
_{m19} &  &  \\ 
&  &  & \ddots  & \ddots  & \ddots  & \ddots  & \ddots  & \ddots  &  \\ 
&  &  &  & \nu _{m7} & \nu _{m8} & \nu _{m9} & \nu _{m10} & \nu _{m7} & \nu
_{m11} \\ 
&  &  &  & \nu _{m18} & \nu _{m19} & \nu _{m20} & \nu _{m21} & \nu _{m22} & 
\nu _{m19}%
\end{bmatrix}%
\end{equation*}
and 
$$\mathbf{x}^{n+1}=[\delta _{-1}^{n+1},\phi _{-1}^{n+1},\delta
_{0}^{n+1},\phi _{0}^{n+1},\ldots ,\delta _{n+1}^{n+1},\phi _{n+1}^{n+1}]^T$$
The system (\ref{mat}) contains $2N+2$ equations and $2N+6$ unknowns. The homogenous Neumann boundary conditions at both ends $\left . \dfrac{\partial V(x,t)}{\partial x}\right |_{x=a}=\left . \dfrac{\partial V(x,t)}{\partial x}\right |_{x=b}=\left . \dfrac{\partial U(x,t)}{\partial x}\right |_{x=a}=\left . \dfrac{\partial U(x,t)}{\partial x}\right |_{x=b}$ are adapted as;
\begin{equation*}
\delta _{-1}=\delta _{1},\text{ }\phi _{-1}=\phi _{1},\text{ }\delta
_{N+1}=\delta _{N-1},\text{ }\phi _{N+1}=\phi _{N-1}
\end{equation*}
in order to convert (\ref{mat}) to a solvable system. The resultant linear equation system is solved using Thomas algorithm.
\section{The Initial State}
The start vector $\mathbf{x^0}$ should be obtained in order to be able to start the iteration in (\ref{mat}). Rearranging the initial and boundary conditions \scriptsize
\begin{equation}
\begin{aligned}
U^{0}(x_0,0)&=\dfrac{\sinh{(\zeta h)}-\zeta h}{2(\zeta h\cosh{(\zeta h)}-\sinh{(\zeta h)})}\delta _{-1}^{0}+\delta _{0}^{0}+\dfrac{\sinh{(\zeta h)}-\zeta h}{2(\zeta h\cosh{(\zeta h)}-\sinh{(\zeta h)})}\delta _{1}^{0} \\ 
U^{0}(x_{m},0)&=\dfrac{\sinh{(\zeta h)}-\zeta h}{2(\zeta h\cosh{(\zeta h)}-\sinh{(\zeta h)})}\delta _{m-1}^{0}+\delta _{m}^{0}+\dfrac{\sinh{(\zeta h)}-\zeta h}{2(\zeta h\cosh{(\zeta h)}-\sinh{(\zeta h)})}\delta _{m+1}^{0},\text{ \  \ }  m=1,2,\ldots ,N-1 \\ 
U^{0}(x_N,0)&=U_{N}^{0}=\dfrac{\sinh{(\zeta h)}-\zeta h}{2(\zeta h\cosh{(\zeta h)}-\sinh{(\zeta h)})}\delta _{N-1}^{0}+\delta _{N}^{0}+\dfrac{\sinh{(\zeta h)}-\zeta h}{2(\zeta h\cosh{(\zeta h)}-\sinh{(\zeta h)})}\delta _{N+1}^{0}
\end{aligned}
\label{B1}
\end{equation}
\normalsize
and
\scriptsize
\begin{equation}
\begin{aligned}
V^{0}(x_0,0)&=\dfrac{\sinh{(\zeta h)}-\zeta h}{2(\zeta h\cosh{(\zeta h)}-\sinh{(\zeta h)})}\phi _{-1}^{0}+\phi _{0}^{0}+\dfrac{\sinh{(\zeta h)}-\zeta h}{2(\zeta h\cosh{(\zeta h)}-\sinh{(\zeta h)})}\phi _{1}^{0} \\ 
V^{0}(x_{m},0)&=\dfrac{\sinh{(\zeta h)}-\zeta h}{2(\zeta h\cosh{(\zeta h)}-\sinh{(\zeta h)})}\phi _{m-1}^{0}+\phi _{m}^{0}+\dfrac{\sinh{(\zeta h)}-\zeta h}{2(\zeta h\cosh{(\zeta h)}-\sinh{(\zeta h)})}\phi _{m+1}^{0},\text{ \  \ }m=1,2,\ldots ,N-1 \\ 
V^{0}(x_N,0)&=\dfrac{\sinh{(\zeta h)}-\zeta h}{2(\zeta h\cosh{(\zeta h)}-\sinh{(\zeta h)})}\phi _{N-1}^{0}+\phi _{N}^{0}+\dfrac{\sinh{(\zeta h)}-\zeta h}{2(\zeta h\cosh{(\zeta h)}-\sinh{(\zeta h)})}\phi _{N+1}^{0}
\end{aligned}
\label{B2}
\end{equation}
\normalsize
gives $N+1$ equations with $N+3$ unknowns. Using 
\begin{equation*}
\begin{aligned}
\delta _{-1}^{0}&=\delta _{1}^{0}+\frac{2(\zeta h\cosh{(\zeta h)}-\sinh{(\zeta h)})}{\zeta (1-\cosh{(\zeta h)})}U_{0}^{\prime } \\
\delta _{N+1}^{0}&=\delta _{N-1}^{0}-\frac{2(\zeta h\cosh{(\zeta h)}-\sinh{(\zeta h)})}{\zeta (1-\cosh{(\zeta h)})}U_{N}^{\prime } \\
\phi _{-1}^{0}&=\phi _{1}^{0}+\frac{2(\zeta \cosh{(\zeta h)}-\sinh{(\zeta h)})}{\zeta (1-\cosh{(\zeta h)})}V_{0}^{\prime } \\
\phi _{N+1}^{0}&=\phi _{N-1}^{0}-\frac{2(\zeta \cosh{(\zeta h)}-\sinh{(\zeta h)})}{\zeta (1-\cosh{(\zeta h)})}V_{N}^{\prime }
\end{aligned}
\end{equation*}
the parameters $\delta _{-1}^{0},\delta _{N+1}^{0}, \phi _{-1}^{0}$ and $\phi _{N+1}^{0}$ can be eliminated from the system (\ref{B1})-(\ref{B2}). The resultant linear equation system;
\begin{equation*}
\left[ 
\begin{tabular}{ccccccc}
$1\smallskip $ & $\tfrac{\tilde{s}-\zeta h}{\zeta h\tilde{c}-\tilde{s}}$ &  &  &  &  &  \\ 
$\tfrac{\tilde{s}-\zeta h}{2(\zeta h\tilde{c}-\tilde{s})}$ & $1$ & $\tfrac{\tilde{s}-\zeta h}{2(\zeta h\tilde{c}-\tilde{s})}$ &  &  &  &  \\ 
&  &  & $\ddots $ &  &  &  \\ 
&  &  &  & $\tfrac{\tilde{s}-\zeta h}{2(\zeta h\tilde{c}-\tilde{s})}\smallskip $ & $1$ & $\tfrac{\tilde{s}-\zeta h}{2(\zeta h\tilde{c}-\tilde{s})}$ \\ 
&  &  &  &  & $\tfrac{\tilde{s}-\zeta h}{\zeta h\tilde{c}-\tilde{s}}$ & $1$%
\end{tabular}
\  \right] \left[ 
\begin{tabular}{c}
$\delta _{0}^{0}\smallskip $ \\ 
$\delta _{1}^{0}$ \\ 
$\vdots $ \\ 
$\delta _{N-1}^{0}\smallskip $ \\ 
$\delta _{N}^{0}$%
\end{tabular}%
\  \right] =\left[ 
\begin{tabular}{c}
$U_{0}-\tfrac{\tilde{s}-\zeta h}{\zeta (1-\tilde{c})}U_{0}^{\prime }\smallskip $ \\ 
$U_{1}^{\prime }$ \\ 
$\vdots $ \\ 
$U_{N-1}^{\prime }\smallskip $ \\ 
$U_{N}-\tfrac{\tilde{s}-\zeta h}{\zeta (\tilde{c}-1)}U_{N}^{\prime }$%
\end{tabular}%
\right] 
\end{equation*}%
and%
\begin{equation*}
\left[ 
\begin{tabular}{ccccccc}
$\smallskip 1$ & $\tfrac{\tilde{s}-\zeta h}{\zeta h\tilde{c}-\tilde{s}}$ &  &  &  &  &  \\ 
$\tfrac{\tilde{s}-\zeta h}{2(\zeta h\tilde{c}-\tilde{s})}$ & $1$ & $\tfrac{\tilde{s}-\zeta h}{2(\zeta h\tilde{c}-\tilde{s})}$ &  &  &  &  \\ 
&  &  & $\ddots $ &  &  &  \\ 
&  &  &  & $\tfrac{\tilde{s}-\zeta }{2(\zeta h\tilde{c}-\tilde{s})}\smallskip $ & $1$ & $\tfrac{\tilde{s}-\zeta h}{2(\zeta h\tilde{c}-\tilde{s})%
}$ \\ 
&  &  &  &  & $\tfrac{\tilde{s}-\zeta h}{\zeta h\tilde{c}-\tilde{s}}$ & $1$%
\end{tabular}%
\  \right] \left[ 
\begin{tabular}{c}
$\phi _{0}^{0}\smallskip $ \\ 
$\phi _{1}^{0}$ \\ 
$\vdots $ \\ 
$\phi _{N-1}^{0}\smallskip $ \\ 
$\phi _{N}^{0}$%
\end{tabular}%
\  \right] =\left[ 
\begin{tabular}{c}
$V_{0}-\tfrac{\tilde{s}-\zeta h}{\zeta (1-\tilde{c})}V_{0}^{\prime }\smallskip $ \\ 
$V_{1}^{\prime }$ \\ 
$\vdots $ \\ 
$V_{N-1}^{\prime }\smallskip $ \\ 
$V_{N}-\tfrac{\tilde{s}-\zeta h}{\zeta (\tilde{c}-1)}V_{N}^{\prime }$%
\end{tabular}%
\right] 
\end{equation*}%
where $\tilde{c}=\cosh{(\zeta h)}$ and $\tilde{s}=\sinh{(\zeta h)}$ is solved by Thomas algorithm for 3-banded systems to obtain the initial vector $\mathbf{x^0}$.
\section{Numerical Examples}

The accuracy and efficiency of the proposed method will be validated by solving initial boundary value problems for RBS and CBS. The numerical results will be compared with the analytical solutions and the Aksoy' s numerical solutions obtained by extended cubic B-spline collocation method (EXCC) with free parameter $\lambda$\cite{aksoy}. The error between the numerical solution $\psi^{num}$ and the analytical solution $\psi^{ana}$ at a time $t$ will be measured by using discrete maximum error norm $L_{\infty }(\psi)$ defined as;
\begin{equation*}
\begin{array}{c}
L_{\infty }(\psi)=\left \vert \psi^{ana}-\psi^{num}\right \vert _{\infty }=\max \limits_{m}\left
\vert \psi^{ana}_{m}-\psi^{num}_{m}\right \vert \\ 
\end{array}%
\end{equation*}
\subsection{Regularized Boussinesq System}
\begin{figure}[ht]
	\centering
		\includegraphics{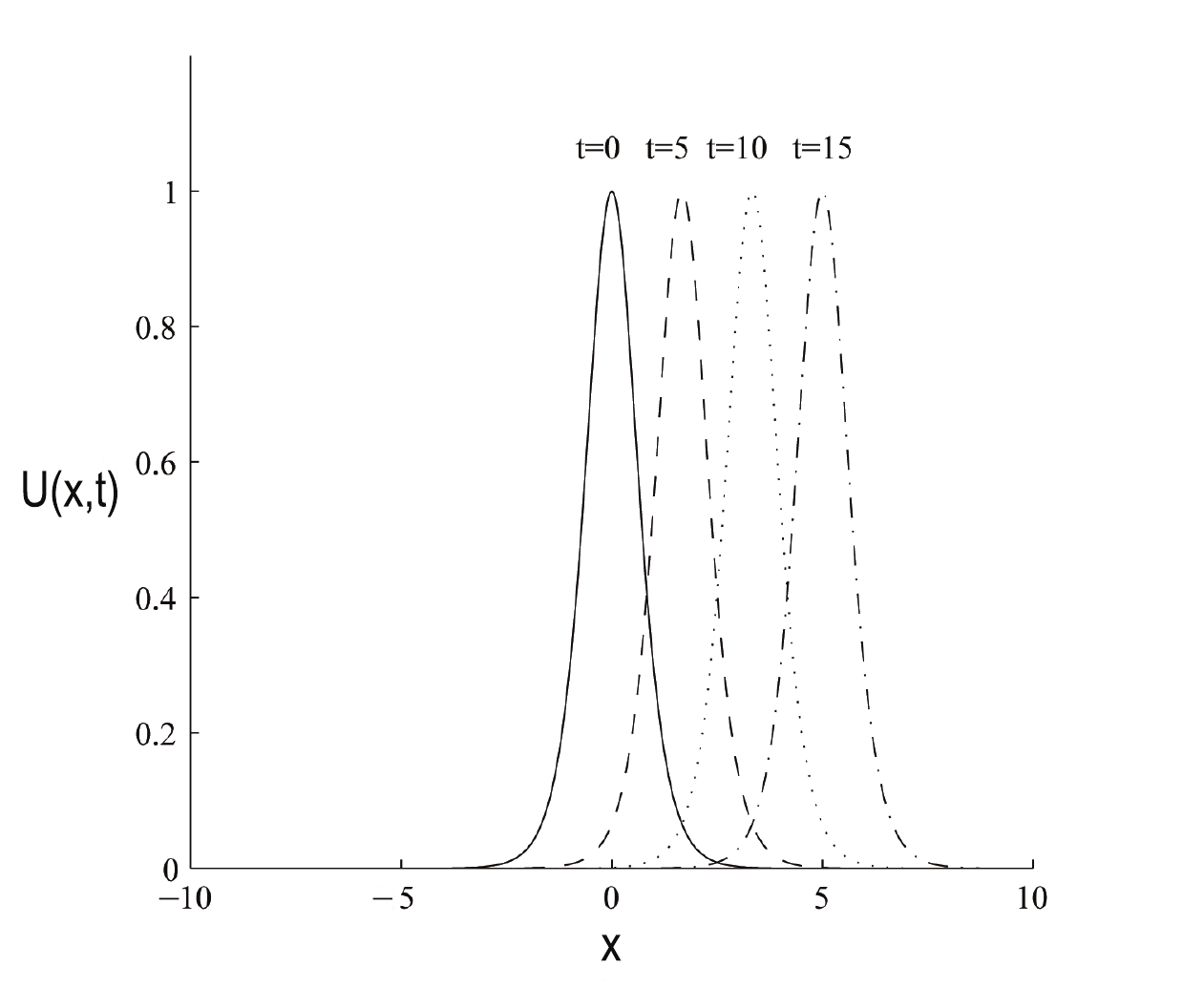}
	\caption{Simulation of wave motion along the horizontal axis}
	\label{fig:2c_tif}
\end{figure}

\begin{figure}[ht]
    \subfigure[Error distribution for $U(x,t)$ at $t=5$]{
   \includegraphics[scale =0.7] {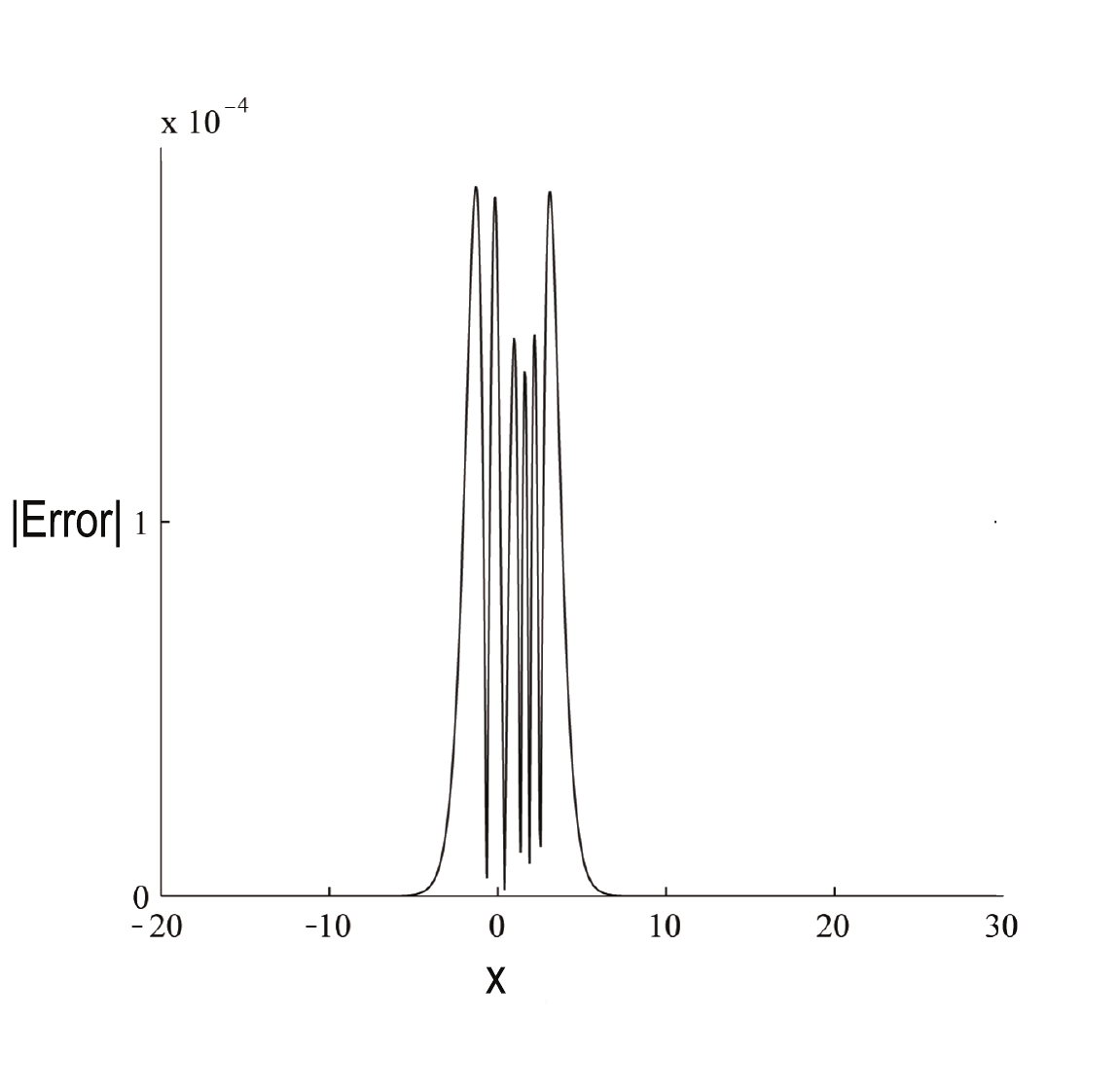}
   \label{fig:2a_eps}
 }
 \subfigure[Error distribution for $V(x,t)$ at $t=5$]{
   \includegraphics[scale =0.7] {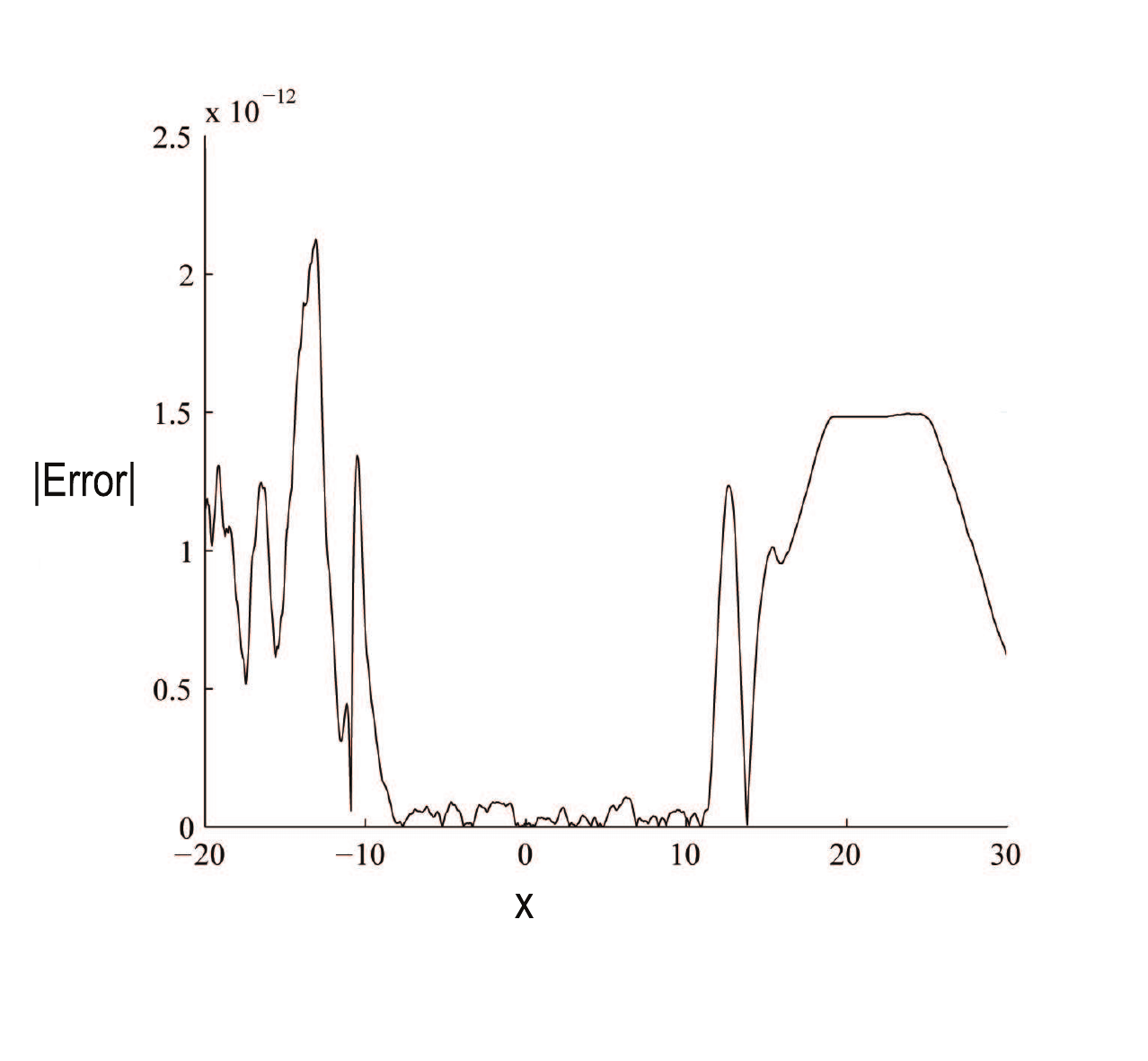}
   \label{fig:2b_eps}
 }
 \caption{Error distribution for both $U(x,t)$ and $V(x,t)$ with  $\zeta =0.0000058339$ and $\Delta t=0.005$ }
\end{figure}

The analytical solution for the RBS is of the form;
\begin{equation*}
\begin{aligned}
V(x,t)&=-1 \\ 
U(x,t)&=(1-\dfrac{\rho}{6})C_s+\dfrac{C_s\rho}{2}\sech^{2}(\dfrac{\sqrt{\rho}}{2}(x+x_{0}-C_st))
\end{aligned}
\end{equation*}
where $x_0$ and $C_s$ real, $\rho$ nonnegative constant\cite{Chena}. This solution represents a traveling wave, the peak of which is located at $x_0$ initially, moving along the horizontal axis with velocity $C_s$. The initial-boundary value problem for the RBS is constructed by combining the initial conditions
\begin{equation}
\begin{aligned}
V(x,0)&=-1 \\ 
U(x,0)&=\sech^{2}(\dfrac{\sqrt{6}}{2}x)
\end{aligned}
\label{RB1}
\end{equation}
and the homogenous Neumann boundary conditions
\begin{equation*}
\begin{aligned}
\left . \dfrac{\partial U(x,t)}{\partial x}\right |_{x=-20}&=0, \, &\left . \dfrac{\partial U(x,t)}{\partial x}\right |_{x=30}=0 \\ 
\left . \dfrac{\partial V(x,t)}{\partial x}\right |_{x=-20}&=0, \, &\left . \dfrac{\partial V(x,t)}{\partial x}\right |_{x=30}=0 
\end{aligned}
\end{equation*}
over the finite interval $[-20,30]$. The initial condition representing a pulse of amplitude $1$ is derived from the analytical solution by a particular choice of $x_{0}=0$, $\rho=6$ and $C_s=1/3$ at the initial time $t=0$. The numerical simulation of wave motion along the horizontal axis is accomplished with the parameters $N=1000$, $\Delta t=0.005$ and $\zeta =0.0000058339$ up to the terminating time $15$, Fig \ref{fig:2c_tif}. The peak of the initial pulse moves from its initial location $x_0=0$ with a constant velocity $C_s=1/3$. At the terminating time $t=15$, the peak reaches $x=5$ without changing its shape. 
For the sake of comparison with Aksoy's work, the designed algorithms are run with fixed number grids $N=1000$, the various parameters $\zeta $ and and time step lengths. The measured maximum errors for both $U(x,t)$ and $V(x,t)$ at the terminating time $t=5$ are depicted in Fig \ref{fig:2a_eps} and Fig \ref{fig:2b_eps}. Even though the maximum error for $U(x,5)$ is located about the peak, the maximum error for $V(x,5)$ is about $x=-15$. 
\noindent
When ECC is compared with EXCC, one see that ECC generates a little bit better numerical results for both various choices of free parameters. When $\Delta t=0.005$, the maximum error is approximately $2\times 10^{-11}$ for $\lambda =0$ and  $3\times 10^{-11}$ for $\lambda =-0.00297$ for the function $V(x,t)$ in EXCC. Usage of the same parameters except $\lambda$, the maximum error is measured less than $10^{-11}$ for various choices of $\zeta $ in ECC. For the function $U(x,t)$, the accuracy of the numerical results are in the same digits for both methods for various values of free parameters $\lambda$ and $\zeta $.   
\begin{table}[H]
  \begin{center}
    \caption{Comparison with Aksoy's work ($L_{\infty}(V)\times 10^5$ at the terminating time $t=5$)}
    \label{tab:table1}
    \begin{tabular}{l|l|l|l|l}
    		\hline\hline
      	$\Delta t$ & EXCC\cite{aksoy} ($\lambda =0$)& ECC ($\zeta =1$)&EXCC\cite{aksoy}& ECC \\ \hline \hline
				$0.5$&$0.000000$&$0.000000$ & $0.000000(\lambda=-0.0642)$ & $0.000000(\zeta =0.0000018762)$ \\ 
				$0.05$ &$0.000000$ & $0.000000$ & $0.000000(\lambda=-0.003627)$ & $0.000000(\zeta =0.0000060571)$ \\ 
				$0.005$ &$0.000002$ &$0.000000$ & $0.000003(\lambda=-0.00297)$ & $0.000000(\zeta =0.0000058339)$ \\ \hline \hline
    \end{tabular}
  \end{center}
\end{table}
\begin{table}[H]
  \begin{center}
    \caption{Comparison with Aksoy's work ($L_{\infty}(U)\times 10^3$ at the terminating time $t=5$)}
    \label{tab:table2}
    \begin{tabular}{l|l|l|l|l}
    		\hline\hline 
      	$\Delta t$ & EXCC\cite{aksoy} ($\lambda =0$)& ECC ($\zeta =1$)&EXCC\cite{aksoy}& ECC \\ \hline\hline
				$0.5$&$23.364410$&$23.890978$ & $4.839046(\lambda =-0.0642)$ & $2.994220(\zeta =0.0000018762)$ \\ 
				$0.05$ &$1.435993$ & $1.554225$ & $0.439634(\lambda=-0.003627)$ & $0.186862(\zeta =0.0000060571)$ \\ 
				$0.005$ &$1.231550$ &$1.351017$ & $0.369655((\lambda=-0.00297)$ & $0.189722(\zeta =0.0000058339)$ \\ \hline \hline
    \end{tabular}
  \end{center}
\end{table}
\subsection{Classical Boussinesq System (CBS)}
\begin{figure}[ht]
	\centering
		\includegraphics{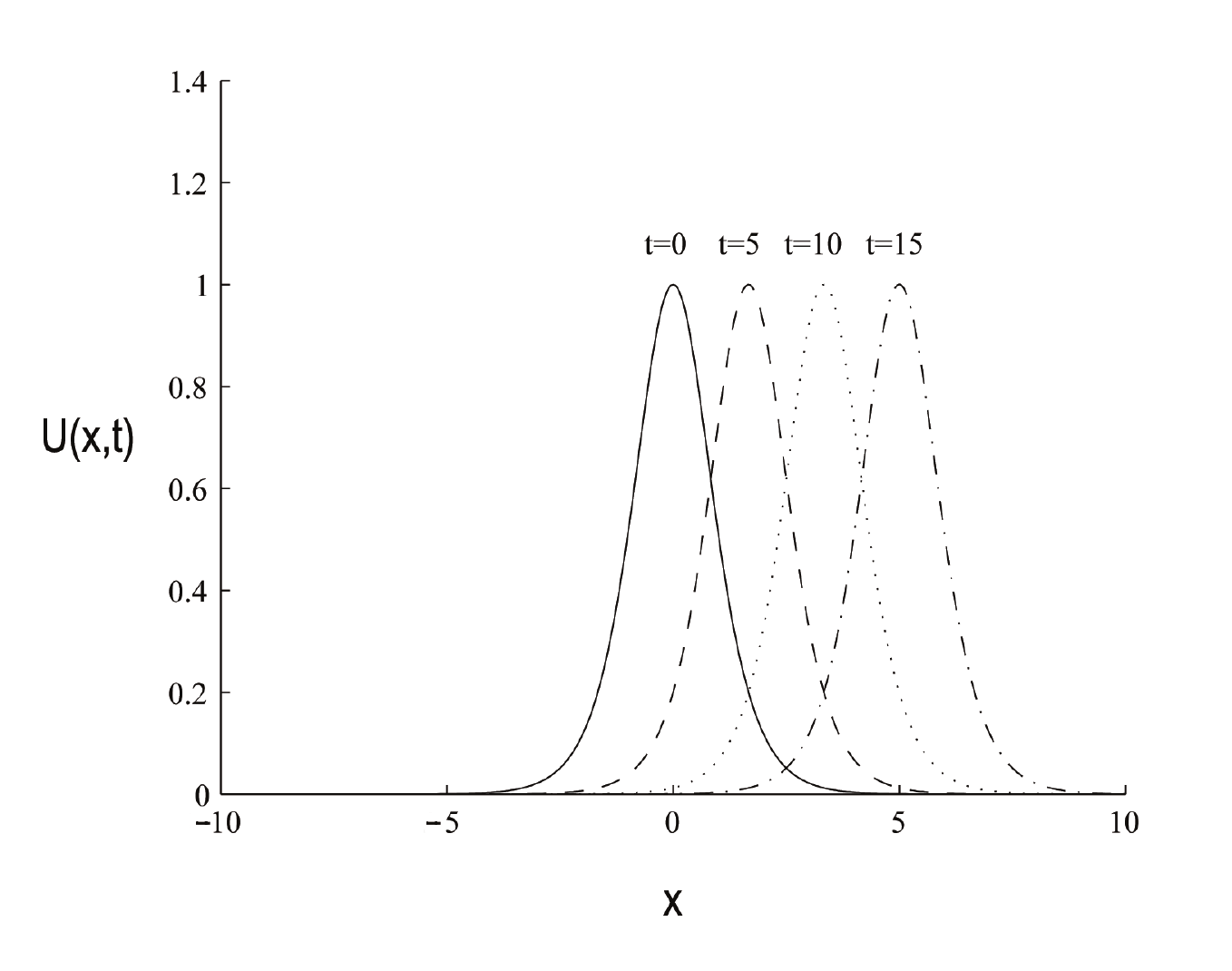}
	\caption{Simulation of traveling wave along the horizontal axis}
	\label{fig:fig3}
\end{figure}

\begin{figure}[ht]
    \subfigure[Error distribution for $U(x,t)$ at $t=5$]{
   \includegraphics[scale =0.7] {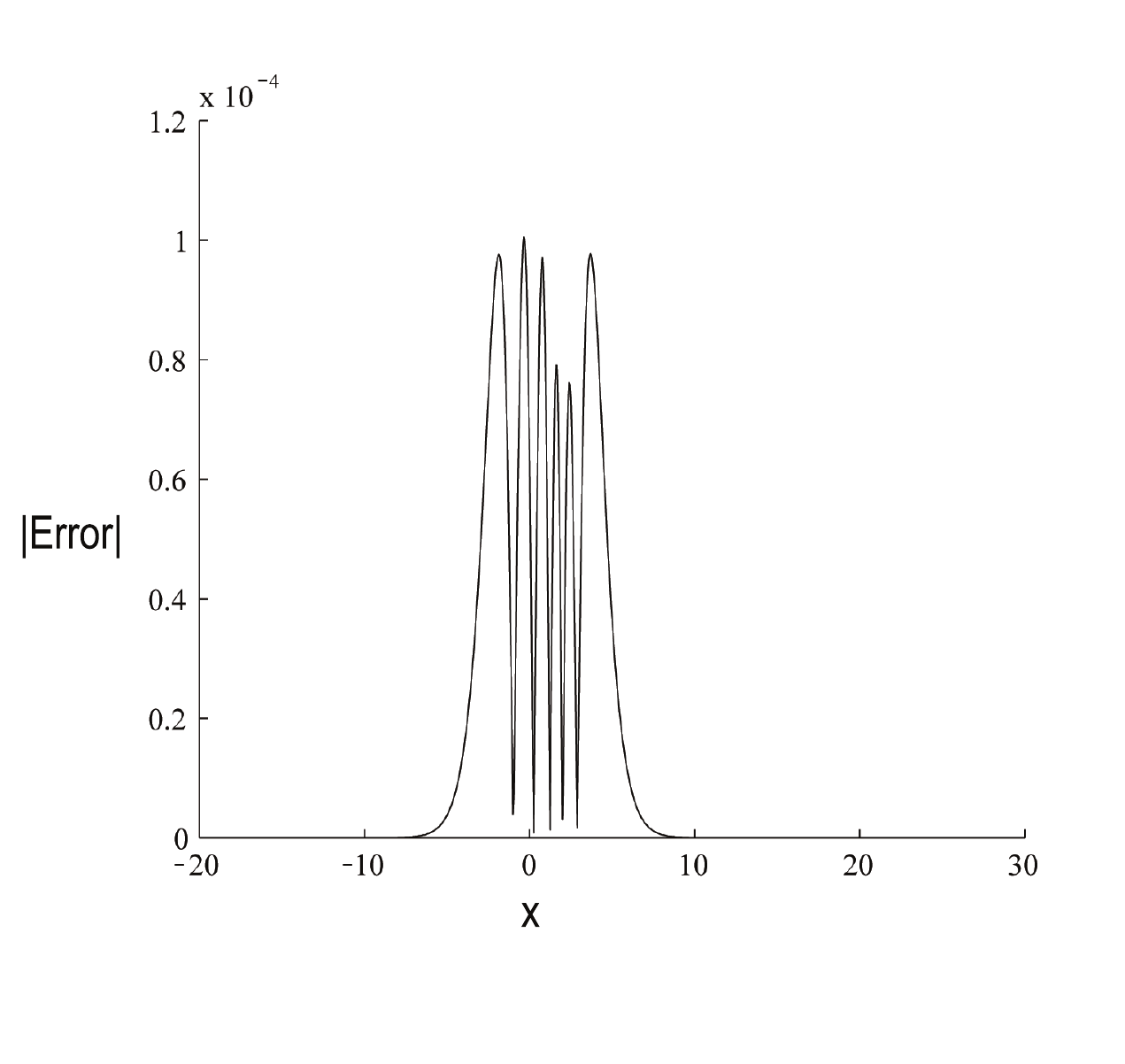}
   \label{fig:4a_eps}
 }
 \subfigure[Error distribution for $V(x,t)$ at $t=5$]{
   \includegraphics[scale =0.7] {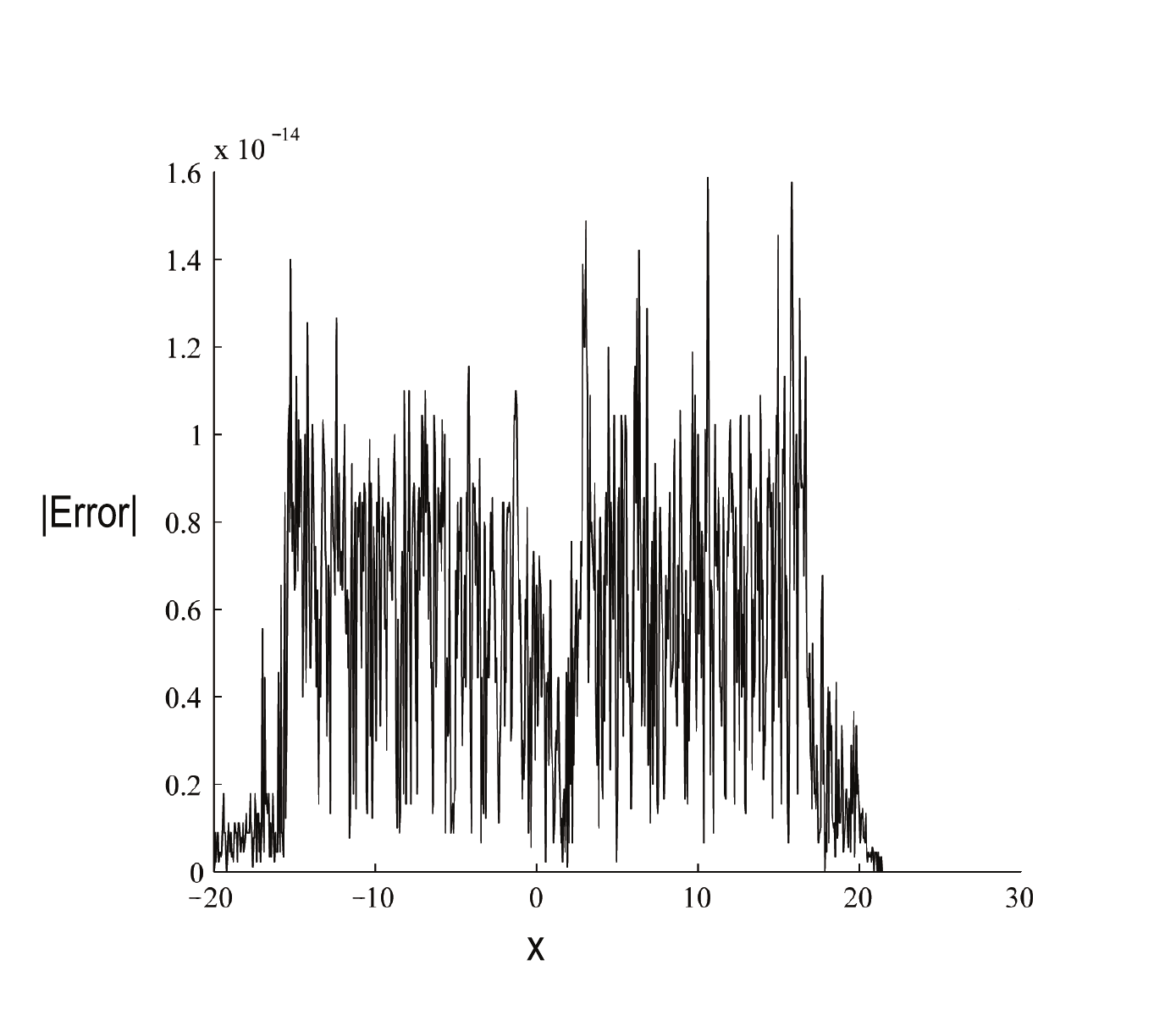}
   \label{fig:4b_eps}
 }
 \caption{Error distribution for both $U(x,t)$ and $V(x,t)$ with  $\zeta =0.0000086530$ and $\Delta t=0.005$ }
\end{figure}

The traveling wave solution for CBS is, 
\begin{equation*}
\begin{aligned}
V(x,t)&=-1 \\ 
U(x,t)&=(1-\dfrac{\rho}{3})C_s+\dfrac{C_s\rho}{2}\sech ^{2}(\dfrac{\sqrt{\rho}}{2}(x+x_{0}-C_st))
\end{aligned}
\end{equation*}
where $x_{0}$ and $C_s$ are real, $\rho$ is non negative constants\cite{Chena}. This solution represents a traveling wave along the horizontal axis with velocity $C_s$ as time goes. 
The initial boundary value problem for CBS is constructed by combining the initial condition
\begin{equation*}
\begin{aligned}
V(x,0)&=-1 \\ 
U(x,0)&=\sech ^{2}(\dfrac{\sqrt{3}}{2}x)
\end{aligned}
\end{equation*}
which is generated by a particular choice of the constants as $x_{0}=0$, $\rho =3$ and $C_s=1/3$ at the initial time $t=0$ in the analytical solution, and the homogenous Neumann boundary conditions covering only first order derivatives at both ends. By those selection of parameters, the analytical solution representing a traveling wave of amplitude $1$ and of initial peak position $0$ moves along the horizontal axis with velocity $1/3$ as time goes. The simulation of the numerical solution is depicted with $N=1000$, $\Delta t=0.005$ and $\zeta =0.0000086530$ up to the terminating time $t=15$, Fig \ref{fig:fig3}. The initial pulse moves along the horizontal axis without changing its shape during the simulation. The peak position becomes $x=t\times C_s=15\times \frac{1}{3}=5$ at the terminating time $t=15$. The error distribution for both function are graphed  with the parameters $N=1000$, $\Delta t=0.005$ and $\zeta =0.0000086530$ at the time $t=5$ in Fig \ref{fig:4a_eps}-\ref{fig:4b_eps}. The maximum errors are measured about $1 \times 10^{-4}$  for $U(x,5)$ and $1 \times 10^{-14}$  for $V(x,5)$.

\noindent
The results are compared with the analytical solutions and Aksoy' s solutions by measuring the error between numerical and analytical solutions for various values of the parameter $\zeta $, Table \ref{tab:table3}- Table \ref{tab:table4}. Even though the accuracy of $V(x,t)$  are almost the same in both methods, the accuracy of $U(x,t)$ appears a little bit better for various values of free parameters in ECC.
\begin{table}[H]
  \begin{center}
    \caption{Comparison with Aksoy's work ($L_{\infty}(V)\times 10^5$ at the terminating time $t=5$)}
    \label{tab:table3}
    \begin{tabular}{l|l|l|l|l}
    		\hline\hline
      	$\Delta t$ & EXCC\cite{aksoy} ($\lambda =0$)& ECC ($\zeta =1$)&EXCC\cite{aksoy}& ECC \\ \hline \hline
				$0.5$&$0.000000$&$0.000000$ & $0.000000(\lambda=-0.0332)$ & $0.000000(\zeta =0.0000027881)$ \\ 
				$0.05$ &$0.000000$ & $0.000000$ & $0.000000(\lambda=-0.00197)$ & $0.000000(\zeta =0.0000080030)$ \\ 
				$0.005$ &$0.000006$ &$0.000006$ & $0.000000(\lambda=-0.0017)$ & $0.000000(\zeta =0.0000086530)$ \\ \hline \hline
    \end{tabular}
  \end{center}
\end{table}
\begin{table}[H]
  \begin{center}
    \caption{Comparison with Aksoy's work ($L_{\infty}(U)\times 10^3$ at the terminating time $t=5$)}
    \label{tab:table4}
    \begin{tabular}{l|l|l|l|l}
    		\hline\hline 
      	$\Delta t$ & EXCC\cite{aksoy} ($\lambda =0$)& ECC ($\zeta =1$)&EXCC\cite{aksoy}& ECC \\ \hline\hline
				$0.5$&$9.915176$&$8.574594$ & $2.279022(\lambda =-0.0339)$ & $0.199469(\zeta =0.0000027881)$ \\ 
				$0.05$ &$0.631048$ & $0.715209$ & $0.206695(\lambda=-0.00197)$ & $0.097838(\zeta =0.0000080030)$ \\ 
				$0.005$ &$0.551502$ &$0.643377$ & $0.182177((\lambda=-0.0017)$ & $0.100474(\zeta =0.0000086530)$ \\ \hline \hline
    \end{tabular}
  \end{center}
\end{table}
\section{Conclusion}
In the study, the exponential cubic B-spline collocation algorithm is designed to solve two initial boundary value problems modeling traveling wave solutions for Regularized Boussinesq and Classical Boussinesq systems. The equations in systems are discretized in time by Crank-Nicolson and finite difference methods. Then, the resultant equations are discretized in space by collocation method based on exponential B-splines. The simulations are obtained successfully by the proposed method. \\
\noindent
The error of numerical results is measured using discrete maximum error norm and compared with the results obtained by Aksoy's method based on extended cubic B-splines. The comparison shows that the results obtained by the proposed method in this study generate more accurate results than Aksoy's results for some problems due to the appropriate choice of free parameter $\zeta$. For the remaining problems, the results are acceptable and as accurate as Aksoy's results.


\begin{thebibliography}{99}
\bibitem{dougalis} Dougalis, V.A., Duran, A., López-Marcos, M.A., Mitsotakis, D.E., A Numerical Study of the Stability of Solitary Waves of the Bona-Smith Family of Boussinesq Systems, Journal of Nonlinear Science, 17, 569-607, 2007.
\bibitem{behzadi} Behzadi S.S., Yildirim A., Application of Quintic B-Spline Collocation Method for Solving the Coupled-BBM System, Middle-East Journal of Scientific Research, 15, 11, 1478-1486, 2013.
\bibitem{anton1} Antonopoulos, D.C., Dougalis V.A., Mitsotakis D.E., Galerkin Approximations of Periodic Solutions of Boussinesq Systems, Bulletin of the Greek Mathematical Society, 57, 13-30, 2010.
\bibitem{karasozen} Karasözen B., Şimşek G., Energy preserving integration of KdV-KdV systems, TWMS Journal of Applied and Engineering Mathematics, 2, 2, 219-227, 2012.
\bibitem{bona2010} Bona, J.L., Dougalis V.A., Mitsotakis, D.E., Numerical solution of KdV-KdV systems of Boussinesq equations I. The numerical scheme and generalized solitary waves, Mathematics and Computers in Simulation, 74, 214-228, 2007. 
\bibitem{anton2} Antonopoulos, D.C., Dougalis V.A., Mitsotakis D.E., Numerical solution of Boussinesq systems of the Bona-Smith family, Applied Numerical Mathematics, 60, 314-336, 2010.
\bibitem{suarez} Su\'{a}rez P.U., Morales, J.H., Numerical Solutions of Two-Way Propagation of Nonlinear Dispersive Waves Using Radial Basis Functions, International Journal of Partial Differential Equations, Article ID 407387, 1-8, 2014.
\bibitem{bona97}  Bona, J.L., Saut, J.-C., Toland, J.F.,  Boussinesq equations for small-amplitude long wavelength water waves, preprint, 1997.
\bibitem{Chena} Chen, M., Exact traveling-wave solutions to bidirectional wave equations, International Journal of Theoretical Physics,
37, 5, 1547-1567, 1998.
\bibitem{BonaChen} Bona, J.L., Chen, M., A Boussinesq system for two way propagation of nonlinear dispersive waves, Physica D, 116, 191-224, 1998.
\bibitem{Bona2002} Bona, J.L., Chen, M., Saut, J.-C., Boussinesq equations and other systems for small-amplitude long waves in nonlinear dispersive media I: derivation and linear theory, Journal of Nonlinear Science, 12, 283-318, 2002.
\bibitem{BS} Boussinesq, J.V., Th\'{e}orie de l'intumescence liquide appel\'{e}e onde solitaire ou de translation se propageant dans un canal rectangulaire, Comptes Rendus de l'Acad\'{e}mie de Sciences, 72, 755-759, 1871.
\bibitem{Chenb} Chen, M., Exact solutions of various Boussinesq systems, Applied Mathematics Letters, 5, 45-49, 1998.
\bibitem{mccartin} McCartin, B.J., Theory of exponential splines, Journal of Approximation Theory, 66, 1, 1-23, 1991.
\bibitem{moh1} Mohammadi, R., Exponential B-Spline Solution of Convection-Diffusion Equations, Applied Mathematics, 4, 933-944, 2013.
\bibitem{moh2} Mohammadi, R., Exponential B-spline collocation method for numerical solution of the generalized regularized long wave equation, Chin. Phys. B, 24, 5, 050206, 2015.
\bibitem{ozkdv} Ersoy, O., Dag, I., The exponential cubic B-spline algorithm for Korteweg-deVries Equation, Advances in Numerical Analysis, Article ID 367056, 2015.
\bibitem{rubin} Rubin, S. G., Graves, R. A., Cubic spline approximation for problems in fluid mechanics, Nasa TR R-436,Washington DC, 1975.
\bibitem{aksoy} Aksoy, A. M., 2012, Numerical Solutions of Some Partial Differential Equations Using the Taylor Collocation-Extended Cubic B-spline Functions, Eskişehir Osmangazi University, Doctoral Dissertation, Eskişehir, 2012.


%\bibitem{Amick} Amick, C.J., 1984, Regularity and uniqueness of solutions to the Boussinesq system of equations, Journal of Differential Equations 54,231-247.
%\bibitem{Schonbek} Schonbek, M.E., 1981, Existence of solutions for the
%Boussinesq sytem of equations, Journal of Differential Equations 42, 325-352.
%\bibitem{Whitnam} Whitham, G.B., 1974, Linear and nonlinear waves, Wiley,
%New York, 636 p.
%\bibitem{PelDou} Pelloni, B., Dougalis, V. A., 2001, Numerical modelling
%of two-way propagation of non-linear dispersive waves, Mathematics and
%Computers in Simulation 55, 595-606.

%\bibitem{BonaSmith} Bona, J. L., Smith, R., 1976, A model for the two-way
%propagation of water waves in a channel, Mathematical Proceedings of the
%Cambridge Philosophical Society 79, 167-182.

%\bibitem{maccartin} McCartin, B. J. , Theory of Exponential Splines, Journal
%of Approximation Theory, Vol. 66, pp. 1-23, 1991.

%\bibitem{sc3} Sakai, M. ,  Usmani, R. A., A class of simple exponential
%B-splines and their application to numerical solution to singular
%perturbation problems, Numer. Math. Vol. 55, pp. 493-500, 1989.

%\bibitem{sc4} Radunovic, D. , Multiresolution exponential B-splines and
%singularly perturbed boundary problem, Numer. Algor., Vol. 47, pp. 191-210,
%2008.

%\bibitem{expcom} R. Delgado-Gonzalo, P. Thevenaz and M. Unser, Exponential
%splines and minimal-support bases for curve representation, Computer Aided
%Geometric Design 29, 109--128, 2012.

%\bibitem{sch} S. Chandra Sekhara Rao and Mukesh Kumar, Exponential B-spline
%collocation method for self-adjoint singularly perturbed boundary value
%problems, Applied Numerical Mathematicals, Vol.58, pp.1572-1581, 2008.

%\bibitem{exp2012} O. Christensen, P. Massopust, Exponential B-splines and
%the partition of unity property, Adv Comput Math, Vol: 37, 301--318, 2012.

%\bibitem{sc2} R. Mohammadi, Exponential B-spline solution \ of
%Convection-Diffusion \ Equations, Applied Mathematics, Vol. 4, pp. 933-944,
%2013.


\end{thebibliography}
\end{document}